\documentclass[twosided,reqno]{amsart}

\usepackage{amsmath}
\usepackage{amssymb}
\usepackage{amsthm}
\usepackage{url}

\theoremstyle{theorem}
\newtheorem{theorem}{\scshape Theorem }[section]
\newtheorem{lemma}[theorem]{\scshape Lemma}

\theoremstyle{definition}

\newtheorem{remark}{\scshape Remark}

\numberwithin{equation}{section}

\begin{document}

\title{Poly-Cauchy numbers and polynomials with umbral calculus viewpoint}

\author{Dae San Kim}
\address{Department of Mathematics, Sogang University, Seoul 121-742, Republic of Korea.}
\email{dskim@sogang.ac.kr}

\author{Taekyun Kim}
\address{Department of Mathematics, Kwangwon University, Seoul 139-701, Republic of Korea}
\email{tkkim@kw.ac.kr}

\maketitle

\begin{abstract}
In this paper, we give some interesting identities of poly-Cauchy numbers and polynomials arising from umbral calculus.
\end{abstract}

\section{Introduction}

For $k \in \mathbb{Z}$, the polylogarithm factorial function is defined by

\begin{equation}\label{1}
Lif_{k}(t)=\sum_{m=0}^{\infty}\frac{t^{m}}{m!(m+1)^{k}},~~(see~[11,12,13]).
\end{equation}

When $k=1$, $Lif_{1}(t)=\sum_{m=0}^{\infty}\frac{t^{m}}{(m+1)!}=\frac{e^{t}-1}{t}.$\\

The poly-Cauchy polynomials of the first kind $C_{n}^{(k)}(x)$ of index $k$ are defined by the generating function to be

\begin{equation}\label{2}
\frac{Lif_{k}(log(1+t))}{(1+t)^{x}}=\sum_{n=0}^{\infty}C_{n}^{(k)}(x)\frac{t^{n}}{n!},~~(k\in \mathbb{Z}),
~~(see[11]).
\end{equation}

For $\lambda \in \mathbb{C}$ with $\lambda \neq 1$, the $n$-th Frobenius-Euler polynomial of order $r$ is defined by the generating function to be

\begin{equation}\label{3}
\left(\frac{1-\lambda}{e^t-\lambda}\right)^{r}e^{xt}=\sum_{n=0} ^{\infty}H_n ^{(r)} (x|\lambda)\frac{t^n}{n!},~~(r \in \mathbb{Z}), {\text{ (see [1,10])}}.
\end{equation}

As is well known, the $n$-th Bernoulli polynomial of order $r$ is given by

\begin{equation}\label{4}
\left(\frac{t}{e^t-1}\right)^{r}e^{xt}=\sum_{n=0} ^{\infty}B_n ^{(r)} (x)\frac{t^n}{n!},~~(r \in \mathbb{Z}), {\text{ (see [1,15])}}.
\end{equation}

The Stirling numbers of the first kind are defined by the generating function to be

\begin{equation}\label{5}
(\log (1+t))^{m}=m!\sum_{l=m} ^{\infty}S_{1}(l,m)\frac{t^l}{l!}~~(m \in \mathbb{Z}_{\geq 0}).
\end{equation}

Thus, by (\ref{5}), we get

\begin{equation}\label{6}
(x)_{n}=x(x-1)\cdots(x-n+1)=\sum_{l=0}^{n}S_{1}(n,l)x^{l},
\end{equation}

and

\begin{equation}\label{7}
x^{(n)}=x(x+1)\cdots(x+n-1)=\sum_{l=0}^{n}(-1)^{n-l}S_{1}(n,l)x^{l},~~(see~[14]).
\end{equation}

Let $\mathbb{C}$ be the complex number field and let ${\mathcal{F}}$ be the set of all formal power series in the variable $t$ over ${\mathbb{C}}$ with

\begin{equation}\label{8}
{\mathcal{F}}=\left\{ \left.f(t)=\sum_{k=0} ^{\infty} a_k\frac{ t^k}{k!}~\right|~ a_k \in {\mathbb{C}} \right\}.
\end{equation}

Let $\mathbb{P}=\mathbb{C}[x]$ and let $\mathbb{P}^{*}$ be the vector space of all linear functionals on ${\mathbb{P}}$.
$\left< L~ \vert ~p(x)\right>$ denotes the action of the linear functional $L$ on the polynomial $p(x)$. For $f(t) \in \mathcal{F}$ with
$f(t)=\sum_{k=0} ^{\infty} a_k\frac{t^{k}}{k!}$, let us define the linear functional on $\mathbb{P}$ by setting

\begin{equation}\label{9}
\left<f(t)|x^n \right>=a_n,~~(n\geq 0), ~~{\text{(see [2,14,15])}}.
\end{equation}

By (\ref{8}) and (\ref{9}), we easily see that

\begin{equation}\label{10}
\left<t^k | x^n \right>=n! \delta_{n,k},~(n,k \geq 0),
\end{equation}

where $\delta_{n,k}$ is the Kronecker symbol. (see [14,15]).

Let us consider $f_{L}(t)=\sum_{k=0}^{\infty}\frac{\left<L | x^n \right>}{k!}t^k$.  By (\ref{10}), we easily see that $\left<f_{L}(t)|x^n \right>=\left<L | x^n \right>$ and
so as linear functionals $L=f_{L}(t)$.  The map $L \mapsto f_{L}(t)$ is a vector space isomorphism from $\mathbb{P}^{*}$ onto $\mathcal{F}$.  Henceforth, $\mathcal{F}$ denotes
both the algebra of formal power series in $t$ and the vector space of all linear functionals on $\mathbb{P}$, and so an element $f(t)$ of $\mathcal{F}$ will be thought of as both a formal power series and a linear functional (see[14,15]).  We call $\mathcal{F}$ the umbral algebra.  The umbral calculus is the study of umbral algebra (see [14]).  The order $o(f(t))$ of a power series
$f(t)(\neq 0)$ is the smallest integer $k$ for which the coefficient of $t^k$ does not vanish.  If $o(f(t))=1$, then $f(t)$ is called a delta series; if $o(g(t))=0$,
then $g(t)$ is called an invertible seires.  Let $f(t),g(t) \in \mathcal{F}$ with $o(f(t))=1$ and $o(g(t))=0$.  Then there exists a unique sequence $S_{n}(x)~(\deg S_{n}(x)=n)$ such that $\left<g(t)f(t)^k | S_n(x) \right>=n! \delta_{n,k}~for~n,k\geq0$. The sequence $S_{n}(x)$ is called the Sheffer sequence for  $(g(t),f(t))$ which is
denoted by $S_n(x)\sim (g(t),f(t)),~(see~[14])$.

For $f(t),g(t) \in \mathcal{F}$ and $p(x) \in \mathbb{P}$, we have

\begin{equation}\label{11}
\left<f(t)g(t) | p(x) \right>=\left<f(t) | g(t)p(x) \right>= \left<g(t) | f(t)p(x) \right>,
\end{equation}

and

\begin{equation}\label{12}
f(t)=\sum_{k=0} ^{\infty} \left<f(t)|x^k\right>\frac{t^k}{k!},~p(x)=\sum_{k=0} ^{\infty}\left<t^k|p(x)\right> \frac{x^k}{k!}.
\end{equation}

By (\ref{12}), we get

\begin{equation}\label{13}
p^{(k)}(0)= \left<t^k|p(x)\right> =\left<1|p^{(k)}(x)\right>~(k \geq 0).
\end{equation}

Thus, by (\ref{13}), we have

\begin{equation}\label{14}
t^{k}p(x)= p^{(k)}(x)= \frac{d^kp(x)}{dx^k},~(see~[14]).
\end{equation}

Let $S_n(x) \sim (g(t), f(t))$.  Then we see that

\begin{equation}\label{15}
\frac{1}{g(\bar{f}(t))}e^{x\bar{f}(t)}= \sum_{n=0}^{\infty}S_n(x)\frac{t^n}{n!},~for~all~x \in \mathbb{C},
\end{equation}

where $\bar{f}(t)$ is the compositional inverse of $f(t)$ with $\bar{f}(f(t))=t$,

\begin{equation}\label{16}
S_{n}(x)=\sum_{j=0}^{n}\frac{1}{j!} \langle \frac{\bar{f}(t)^{j}}{g(\bar{f}(t))} \vert x^{n} \rangle x^{j},
\end{equation}

and

\begin{equation}\label{17}
f(t)S_n(x)=nS_{n-1}(x),~~(n \geq 0),~(see~[14,15]).
\end{equation}

As is well known, the transfer formula for $p_{n}(x) \sim (1,f(t))$, $q_{n} \sim (1,g(t)),~(n \geq 1)$, is given by

\begin{equation}\label{18}
q_n(x)=x \left(\frac{f(t)}{g(t)}\right)^{n}x^{-1}p_{n}(x),~~(see~[14]).
\end{equation}

Let $S_n(x) \sim (g(t),f(t))$. Then it is known that

\begin{equation}\label{19}
S_{n+1}(x)=\left(x-\frac{g'(t)}{g(t)}\right)\frac{1}{f'(t)}S_n(x),~~(see~[14]).
\end{equation}

For $S_n(x) \sim (g(t),f(t))$, $r_n(x) \sim (h(t),l(t))$ we have

\begin{equation}\label{20}
S_{n}(x)=\sum_{m=0}^{n}C_{n,m}r_n(x),
\end{equation}

where

\begin{equation}\label{21}
C_{n,m}=\frac{1}{m!} \langle \frac{h(\bar{f}(t))}{g(\bar{f}(t))}l(\bar{f}(t))^m \vert x^n \rangle, ~~(see~[14]).
\end{equation}

Finally, we note that $e^{yt}p(x)=p(x+y),~(p(x) \in \mathbb{P})$.

In this paper, we investigate some properties of poly-Cauchy numbers and polynomials with umbral calculus viewpoint.  From our investigation, we derive some interesting identities of poly-Cauchy numbers and polynomials.

\section{Poly-Cauchy numbers and polynomials}

From (\ref{2}), we note that $C_{n}^{(k)}(x)$ is the Sheffer sequence for the pair $(g(t)=\frac{1}{Lif_k(-t)},~f(t)=e^{-t}-1)$, that is,

\begin{equation}\label{22}
C_{n}^{(k)}(x) \sim \left(\frac{1}{Lif_k(-t)},e^{-t}-1 \right),~(k \in \mathbb{Z},~n \geq 0).
\end{equation}

When $x=0$, $C_{n}^{(k)}=C_{n}^{(k)}(0)$ is called the $n$-th poly-Cauchy number of the first kind with index $k$.

Thus, we see that

\begin{equation}\label{23}
Lif_{k}(log(1+t))=\sum_{m=0}^{\infty}C_{n}^{(k)}\frac{t^{n}}{n!},~~(see~[12,13]).
\end{equation}

By (\ref{16}) and (\ref{22}), we easily get

\begin{equation}\label{24}
C_{n}^{(k)}(x)=\sum_{j=0}^{n}\frac{1}{j!} \langle Lif_{k}(log(1+t))(-log(1+t))^{j} \vert x^{n} \rangle x^{j}.
\end{equation}

Now, we compute.

\begin{equation}\label{25}
\begin{split}
&\langle Lif_{k}(log(1+t))(-log(1+t))^{j} \vert x^{n} \rangle =(-1)^{j}\sum_{m=0}^{\infty}\frac{1}{m!(m+1)^{k}}\langle (log(1+t))^{m+j} \vert x^{n} \rangle\\
&=(-1)^{j}\sum_{m=0}^{n-j}\frac{1}{m!(m+1)^{k}}\sum_{l=0}^{n-j-m}\frac{(m+j)!}{(l+m+j)!}S_{1}
(l+m+j,m+j)(l+m+j)!\delta_{n,l+m+j}\\
&=(-1)^{j}\sum_{m=0}^{n-j}\frac{(m+j)!}{m!(m+1)^{k}}S_{1}(n,m+j).\\
\end{split}
\end{equation}

Thus, by (\ref{24}) and (\ref{25}), we get

\begin{equation}\label{26}
\begin{split}
C_{n}^{(k)}(x)&=\sum_{j=0}^{n}\{(-1)^{j}\sum_{m=0}^{n-j}\frac{\binom{m+j}{m}}{(m+1)^{k}}S_{1}(n,m+j)\}x^{j}\\
&=\sum_{j=0}^{n}\{(-1)^{j}\sum_{m=j}^{n}\frac{\binom{m}{j}}{(m-j+1)^{k}}S_{1}(n,m)\}x^{j}\\
&=\sum_{m=0}^{n}S_{1}(n,m)\sum_{j=0}^{m}\binom{m}{j}\frac{(-x)^{j}}{(m-j+1)^{k}}.\\
\end{split}
\end{equation}

From (\ref{22}), we have

\begin{equation}\label{27}
\frac{1}{Lif_k(-t)}C_{n}^{(k)}(x)\sim (1,e^{-t}-1),~x^n \sim (1,t).
\end{equation}

By (\ref{18}) and (\ref{22}), for $n\geq1$ we get

\begin{equation}\label{28}
\frac{1}{Lif_k(-t)}C_{n}^{(k)}(x) = x \left( \frac{t}{e^{-t}-1}\right)^{n}x^{-1}x^n=(-1)^nx \left( \frac{te^t}{e^{t}-1}\right)^{n} x^{n-1}.
\end{equation}

Now, we observe that

\begin{equation}\label{29}
\begin{split}
\left( \frac{te^t}{e^{t}-1}\right)^{n}=\left( \frac{-t}{e^{-t}-1}\right)^{n}&= \left( \sum_{l_{1}=0}^{\infty}\frac{(-1)^{l_{1}}t^{l_{1}}}{l_{1}!}B_{l_{1}}\right) \times \cdots \times \left( \sum_{l_{n}=0}^{\infty}\frac{(-1)^{l_{n}}B_{l_{n}}}{l_{n}!}t^{l_{n}}\right)\\
&=\sum_{l=0}^{\infty}\left( \sum_{l_{1}+\cdots + l_{n}=l}(-1)^{l}\binom{l}{l_{1},\cdots , l_{n}}B_{l_{1}} \cdots B_{l_{n}}\right) \frac{t^{l}}{l!},\\
\end{split}
\end{equation}

where $B_{n}$ is the $n$-th ordinary Bernoulli number and
 $\binom{n}{l_{1},\cdots , l_{n}}=\frac{n!}{l_{1}!l_{2}! \cdots l_{n}!}$.

From (\ref{28}) and (\ref{29}), we have

\begin{equation}\label{30}
\begin{split}
&C_{n}^{(k)}(x)\\
&=(-1)^{n}\sum_{l=0}^{n-1} \sum_{l_{1}+\cdots + l_{n}=l}(-1)^{l}\binom{n-1}{l_{1},\cdots , l_{n},n-1-l}B_{l_{1}} \cdots B_{l_{n}}Lif_k(-t)x^{n-l}\\
&=(-1)^{n}\sum_{l=0}^{n-1} \sum_{l_{1}+\cdots + l_{n}=l}\sum_{m=0}^{n-l}(-1)^{m+l}\binom{n-1}{l_{1},\cdots , l_{n},n-1-l}\binom{n-l}{m}\\
&\qquad\qquad\qquad\qquad\qquad\qquad\qquad\qquad  \times \frac{1}{(m+1)^k}B_{l_{1}} \cdots B_{l_{n}}x^{n-l-m}\\
&=(-1)^{n}\sum_{l=0}^{n-1} \sum_{l_{1}+\cdots + l_{n}=l}\sum_{j=0}^{n-l}(-1)^{n-j}\binom{n-1}{l_{1},\cdots , l_{n},n-1-l}\binom{n-l}{j}\frac{ B_{l_{1}} \cdots B_{l_{n}}}{(n-l-j+1)^k}x^{j}\\
&=\sum_{l=0}^{n-1} \sum_{l_{1}+\cdots + l_{n}=l}\binom{n-1}{l_{1},\cdots , l_{n},n-1-l}\frac{ B_{l_{1}} \cdots B_{l_{n}}}{(n-l+1)^k}\\
& \qquad +\sum_{j=1}^{n}\left\{\sum_{l=0}^{n-j}\sum_{l_{1}+\cdots + l_{n}=l}(-1)^{j}\binom{n-1}{l_{1},\cdots , l_{n},n-1-l}\binom{n-l}{j}\frac{ B_{l_{1}} \cdots B_{l_{n}}}{(n-l-j+1)^k}\right\}x^{j}.\\
\end{split}
\end{equation}

Therefore, by (\ref{30}), we obtain the following theorem

\begin{theorem}\label{thm1}
For $k \in \mathbb{Z},~n \geq 1$, we have
\begin{equation*}
\begin{split}
C_{n}^{(k)}(x)&=\sum_{l=0}^{n-1} \sum_{l_{1}+\cdots + l_{n}=l}\binom{n-1}{l_{1},\cdots , l_{n},n-1-l}\frac{ B_{l_{1}} \cdots B_{l_{n}}}{(n-l+1)^k}\\
& +\sum_{j=1}^{n}\left\{\sum_{l=0}^{n-j}\sum_{l_{1}+\cdots + l_{n}=l}(-1)^{j}\binom{n-1}{l_{1},\cdots , l_{n},n-1-l}\binom{n-l}{j}\frac{ B_{l_{1}} \cdots B_{l_{n}}}{(n-l-j+1)^k}\right\}x^{j}.\\
\end{split}
\end{equation*}
\end{theorem}

From (\ref{28}), we have

\begin{equation}\label{31}
\begin{split}
\frac{1}{Lif_k(-t)}C_{n}^{(k)}(x) &= (-1)^nx \left( \frac{te^t}{e^{t}-1}\right)^{n} x^{n-1}=(-1)^nx \left( t+ \frac{t}{e^{t}-1}\right)^{n} x^{n-1}\\
&=(-1)^nx \sum_{a=0}^{n}\binom{n}{a}t^{n-a} \left( \frac{t}{e^{t}-1}\right)^{a} x^{n-1}\\
&=(-1)^nx \sum_{a=0}^{n}\binom{n}{a}t^{n-a} B_{n-1}^{(a)}(x)\\
&=(-1)^n \sum_{a=1}^{n}\binom{n}{a}(n-1)_{n-a}x B_{a-1}^{(a)}(x)\\
&=n!(-1)^n \sum_{a=1}^{n}\frac{1}{a!}\binom{n-1}{a-1}x B_{a-1}^{(a)}(x).\\
\end{split}
\end{equation}

Thus, by (\ref{31}), we get

\begin{equation}\label{32}
\begin{split}
C_{n}^{(k)}(x)&=(-1)^n n! \sum_{a=1}^{n}\sum_{l=0}^{a-1}\frac{1}{a!}\binom{n-1}{a-1}\binom{a-1}{l} B_{a-1-l}^{(a)}Lif_k(-t)x^{l+1}\\
&=(-1)^n n! \sum_{a=1}^{n}\sum_{l=0}^{a-1}\sum_{m=0}^{l+1}\frac{1}{a!}\binom{n-1}{a-1}\binom{a-1}{l}\binom{l+1}{m}\frac{(-1)^{l+1-m}}{(l+2-m)^{k}} B_{a-1-l}^{(a)}x^{m}\\
&=(-1)^n n!\left\{ \sum_{a=1}^{n}\sum_{l=0}^{a-1}\frac{(-1)^{l+1}}{a!}\binom{n-1}{a-1}\binom{a-1}{l}\frac{B_{a-1-l}^{(a)}}{(l+2)^{k}}\right.\\
&\left. + \sum_{a=1}^{n}\sum_{l=0}^{a-1}\sum_{m=1}^{l+1}\frac{1}{a!}\binom{n-1}{a-1}\binom{a-1}{l}\binom{l+1}{m}\frac{(-1)^{l+1-m}}{(l+2-m)^{k}} B_{a-1-l}^{(a)}x^{m}\right\}\\
&=(-1)^n n!\left\{ \sum_{a=1}^{n}\sum_{l=0}^{a-1}\frac{(-1)^{l+1}}{a!}\binom{n-1}{a-1}\binom{a-1}{l}\frac{B_{a-1-l}^{(a)}}{(l+2)^{k}}\right.\\
&\left. + \sum_{m=1}^{n}\sum_{a=m}^{n}\sum_{l=m-1}^{a-1}\frac{(-1)^{l+1-m}}{a!}\binom{n-1}{a-1}\binom{a-1}{l}\binom{l+1}{m}\frac{ B_{a-1-l}^{(a)}}{(l+2-m)^{k}}x^{m}\right\}.\\
\end{split}
\end{equation}

Therefore, by (\ref{26}),(\ref{30}) and (\ref{32}), we obtain the following theorem.

\begin{theorem}\label{thm2}
For $n \geq 1,~ 1 \leq j \leq n$, we have
\begin{equation*}
\begin{split}
&\sum_{m=j}^{n}\frac{\binom{m}{j}}{(m-j+1)^{k}}S_{1}(n,m)\\
&=\sum_{l=0}^{n-j}\sum_{l_{1}+\cdots + l_{n}=l}\binom{n-1}{l_{1},\cdots , l_{n},n-1-l}\binom{n-l}{j}\frac{ B_{l_{1}} \cdots B_{l_{n}}}{(n-l-j+1)^k}\\
&=(-1)^n n! \sum_{a=j}^{n}\sum_{l=j-1}^{a-1}\frac{(-1)^{l+1-j}}{a!}\binom{n-1}{a-1}\binom{a-1}{l}\binom{l+1}{j}\frac{B_{a-1-l}^{(a)}}{(l+2-j)^{k}}.\\
\end{split}
\end{equation*}
Moreover,
\begin{equation*}
\begin{split}
C_{n}^{(k)}&=\sum_{m=0}^{n}\frac{S_{1}(n,m)}{(m+1)^{k}}=\sum_{l=0}^{n-1}\sum_{l_{1}+\cdots + l_{n}=l}\binom{n-1}{l_{1},\cdots , l_{n},n-1-l}\frac{ B_{l_{1}} \cdots B_{l_{n}}}{(n-l+1)^k}\\
&=(-1)^n n! \sum_{a=1}^{n}\sum_{l=0}^{a-1}\frac{(-1)^{l+1}}{a!}\binom{n-1}{a-1}\binom{a-1}{l}\frac{B_{a-1-l}^{(a)}}{(l+2)^{k}}.\\
\end{split}
\end{equation*}
\end{theorem}
where $k \in \mathbb{Z},~n \geq 1$.

From (\ref{7}), we note that

\begin{equation}\label{33}
\frac{1}{(1-t)^{x}}=\sum_{n=0}^{\infty}(-x)_{n}\frac{(-t)^{n}}{n!}=\sum_{n=0}^{\infty}x^{(n)}\frac{t^{n}}{n!}.
\end{equation}

By (\ref{15}) and (\ref{33}), we get

\begin{equation}\label{34}
x^{(n)}=\sum_{m=0}^{n}(-1)^{n-m}S_{1}(n,m)x^{m} \sim (1,1-e^{-t}),
\end{equation}

and

\begin{equation}\label{35}
(-1)^{n}x^{(n)}=\sum_{m=0}^{n}(-1)^{m}S_{1}(n,m)x^{m} \sim (1,e^{-t}-1).
\end{equation}

Thus, by (\ref{27}) and (\ref{35}), we get

\begin{equation}\label{36}
\frac{1}{Lif_k(-t)}C_{n}^{(k)}(x)=(-1)^{n}x^{(n)}  \Leftrightarrow C_{n}^{(k)}(x)=(-1)^{n}Lif_k(-t)x^{(n)}.
\end{equation}

From (\ref{36}), we have

\begin{equation}\label{37}
\begin{split}
C_{n}^{(k)}(x)&=Lif_k(-t)(-1)^{n}x^{(n)}=\sum_{m=0}^{n}(-1)^{m}S_{1}(n,m)Lif_k(-t)x^{m}\\
&=\sum_{m=0}^{n}(-1)^{m}S_{1}(n,m) \sum_{a=0}^{\infty}\frac{(-1)^{a}}{a!(a+1)^k}t^{a}x^{m}\\
&=\sum_{m=0}^{n}(-1)^{m}S_{1}(n,m) \sum_{a=0}^{m}\frac{(-1)^{a}}{a!(a+1)^k}(m)_{a}x^{m-a}\\
&=\sum_{m=0}^{n}\sum_{j=0}^{m}(-1)^{m}S_{1}(n,m)\frac{(-1)^{m-j}}{(m-j)!(m-j+1)^k}(m)_{m-j}x^{j}\\
&=\sum_{m=0}^{n}\sum_{j=0}^{m}S_{1}(n,m)\binom{m}{j}\frac{(-x)^{j}}{(m-j+1)^k}.\\
\end{split}
\end{equation}

It is well known that the Sheffer identity is given by

\begin{equation}\label{38}
S_{n}(x+y)=\sum_{j=0}^{n}\binom{n}{j}S_{j}(x)P_{n-j}(y),~~(see~[14]),
\end{equation}
where $S_{n}(x) \sim (g(t),f(t))$ and $P_{n}(x)=g(t)S_{n}(x)$.

By (\ref{38}), we easily get

\begin{equation}\label{39}
C_{n}^{(k)}(x+y)=\sum_{j=0}^{n}(-1)^{n-j}\binom{n}{j}C_{j}^{(k)}(x)y^{(n-j)},
\end{equation}
where $y^{(n)}=y(y+1) \cdots (y+n-1)$.

From  (\ref{19}) and  (\ref{22}), we have

\begin{equation}\label{40}
\begin{split}
C_{n+1}^{(k)}(x)&=\left(x-\frac{Lif'_k(-t)}{Lif_k(-t)}\right)(-e^{t})C_{n}^{(k)}(x)\\
&=e^{t}\frac{Lif'_k(-t)}{Lif_k(-t)}C_{n}^{(k)}(x)-xC_{n}^{(k)}(x+1).\\
\end{split}
\end{equation}

Now, we compute

\begin{equation}\label{41}
\begin{split}
\frac{Lif'_k(-t)}{Lif_k(-t)}C_{n}^{(k)}(x)&=Lif'_k(-t)\left(\frac{1}{Lif_k(-t)}C_{n}^{(k)}(x) \right)\\
&=Lif'_k(-t)(-1)^{n}x^{(n)}=(-1)^{n}Lif'_k(-t)x^{(n)}\\
&=(-1)^{n}\sum_{l=0}^{n}(-1)^{n-l}S_{1}(n,l)Lif'_k(-t)x^{l}.\\
\end{split}
\end{equation}

By the definition of the polylogarithm factorial function, we get

\begin{equation}\label{42}
Lif'_k(-t)=\sum_{n=0}^{\infty}\frac{(-1)^{n}t^{n}}{n!(n+2)^{k}}.
\end{equation}

From (\ref{41}) and  (\ref{42}), we can derive

\begin{equation}\label{43}
\begin{split}
\frac{Lif'_k(-t)}{Lif_k(-t)}C_{n}^{(k)}(x)&=(-1)^{n}\sum_{l=0}^{n}(-1)^{n-l}S_{1}(n,l)\sum_{m=0}^{l}(-1)^{m}\binom{l}{m}\frac{x^{l-m}}{(m+2)^{k}}\\
&=\sum_{l=0}^{n}S_{1}(n,l)\sum_{j=0}^{l}\frac{(-1)^{j}\binom{l}{j}}{(l-j+2)^{k}}x^{j}.\\
\end{split}
\end{equation}

Thus, by (\ref{40}) and  (\ref{43}), we get

\begin{equation}\label{44}
\begin{split}
C_{n+1}^{(k)}(x)&=e^{t}\sum_{l=0}^{n}S_{1}(n,l)\sum_{j=0}^{l}\frac{(-1)^{j}\binom{l}{j}}{(l-j+2)^{k}}x^{j}-xC_{n}^{(k)}(x+1)\\
&=\sum_{l=0}^{n}S_{1}(n,l)\sum_{j=0}^{l}\frac{(-1)^{j}\binom{l}{j}}{(l-j+2)^{k}}(x+1)^{j}-xC_{n}^{(k)}(x+1).\\
\end{split}
\end{equation}

Therefore, we obtain the following theorem.

\begin{theorem}\label{thm3}
For $k \in \mathbb{Z},~n \geq 0$, we have
\begin{equation*}
C_{n+1}^{(k)}(x)=\sum_{l=0}^{n}S_{1}(n,l)\sum_{j=0}^{l}\frac{(-1)^{j}\binom{l}{j}}{(l-j+2)^{k}}(x+1)^{j}-xC_{n}^{(k)}(x+1).
\end{equation*}
\end{theorem}

For $f(t) \in \mathcal{F}$ and $p(x) \in \mathbb{P}$, we note that

\begin{equation}\label{45}
\langle f(t) \vert xp(x) \rangle =\langle \partial_{t} f(t) \vert p(x) \rangle ,~~(see~[14]),
\end{equation}
where $\partial_{t}f(t) = \frac{df(t)}{dt}$.

By (\ref{10}) and  (\ref{45}), we get

\begin{equation}\label{46}
\begin{split}
C_{n}^{(k)}(y)&=\langle \sum_{l=0}^{\infty}C_{l}^{(k)}(y)\frac{t^{l}}{l!} \vert x^{n} \rangle =\langle \frac{Lif_k(\log(1+t))}{(1+t)^{y}} \vert xx^{n-1} \rangle\\
&=\langle \partial_{t}( Lif_k(\log(1+t))(1+t)^{-y}) \vert x^{n-1} \rangle \\
& =\langle( \partial_{t} Lif_k(\log(1+t)))(1+t)^{-y} \vert x^{n-1} \rangle +\langle Lif_k(\log(1+t))\partial_{t}(1+t)^{-y} \vert x^{n-1} \rangle\\
&=\langle Lif'_k(\log(1+t))(1+t)^{-y-1} \vert x^{n-1} \rangle-yC_{n-1}^{(k)}(y+1).\\
\end{split}
\end{equation}

It is easy to show that

\begin{equation}\label{47}
(tLif_{k}(t))'=Lif_{k-1}(t),~(tLif_{k}(t))'=Lif_{k}(t)+tLif'_{k}(t).
\end{equation}

Thus, by (\ref{47}), we get

\begin{equation}\label{48}
Lif'_{k}(t)=\frac{Lif_{k-1}(t)-Lif_{k}(t)}{t}
\end{equation}

From (\ref{48}), we can derive the following equation:

\begin{equation}\label{49}
\begin{split}
& \langle Lif'_k(\log(1+t))(1+t)^{-y-1} \vert x^{n-1} \rangle\\
&= \langle \frac{Lif_{k-1}(\log(1+t))-Lif_{k}(\log(1+t))}{t}(1+t)^{-y-1} \vert \frac{t}{\log(1+t)}x^{n-1} \rangle\\
&=\sum_{l=0}^{n-1}\binom{n-1}{l}B_{l}^{(l)}(1) \langle \frac{Lif_{k-1}(\log(1+t))-Lif_{k}(\log(1+t))}{t}(1+t)^{-y-1} \vert x^{n-1-l} \rangle\\
&=\sum_{l=0}^{n-1}\binom{n-1}{l}B_{l}^{(l)}(1) \langle \frac{Lif_{k-1}(\log(1+t))-Lif_{k}(\log(1+t))}{t}(1+t)^{-y-1} \vert \frac{1}{n-l}tx^{n-l}  \rangle\\
&=\sum_{l=0}^{n-1}\binom{n-1}{l}\frac{B_{l}^{(l)}(1)}{n-l}  \langle \left( Lif_{k-1}(\log(1+t))-Lif_{k}(\log(1+t))\right)(1+t)^{-y-1} \vert x^{n-l} \rangle\\
&=\frac{1}{n}\sum_{l=0}^{n-1}\binom{n}{l}B_{l}^{(l)}(1)\{C_{n-l}^{(k-1)}(y+1)-C_{n-l}^{(k)}(y+1)\}.\\
\end{split}
\end{equation}

Therefore, by (\ref{46}) and  (\ref{49}), we obtain the following theorem.

\begin{theorem}\label{thm4}
For $k \in \mathbb{Z},~n \geq 0$, we have
\begin{equation*}
C_{n}^{(k)}(x)=-xC_{n-1}^{(k)}(x+1)+\frac{1}{n} \sum_{l=0}^{n}\binom{n}{l}B_{l}^{(l)}(1)\{C_{n-l}^{(k-1)}(x+1)-C_{n-l}^{(k)}(x+1)\}.
\end{equation*}
\end{theorem}

For $n \geq m \geq 1$, we evaluate

\begin{equation}\label{50}
\langle (\log(1+t))^{m}Lif_{k}(\log(1+t)) \vert x^{n} \rangle
\end{equation}

in two different ways.

On the one hand, we get

\begin{equation}\label{51}
\begin{split}
&\langle (\log(1+t))^{m}Lif_{k}(\log(1+t)) \vert x^{n} \rangle =\langle Lif_{k}(\log(1+t)) \vert (\log(1+t))^{m} x^{n} \rangle\\
&=\sum_{l=0}^{n-m}\frac{m!}{(l+m)!}S_{1}(l+m,m)(n)_{l+m}\langle Lif_{k}(\log(1+t)) \vert  x^{n-l-m} \rangle\\
&=\sum_{l=0}^{n-m}m!\binom{n}{l+m}S_{1}(l+m,m)C_{n-l-m}^{(k)}.\\
\end{split}
\end{equation}

On the other hand, we have

\begin{equation}\label{52}
\begin{split}
&\langle (\log(1+t))^{m}Lif_{k}(\log(1+t)) \vert x^{n} \rangle =\langle (\log(1+t))^{m}Lif_{k}(\log(1+t)) \vert xx^{n-1} \rangle\\
&=\langle \partial_{t}(\log(1+t))^{m}Lif_{k}(\log(1+t)) \vert x^{n-1} \rangle.\\
\end{split}
\end{equation}

Now, we observe that

\begin{equation}\label{53}
\begin{split}
& \partial_{t}((\log(1+t))^{m}Lif_{k}(\log(1+t)) )=\partial_{t}\{(\log(1+t))^{m-1}\log(1+t)Lif_{k}(\log(1+t))\}\\
&=(\partial_{t}(\log(1+t))^{m-1})\log(1+t)Lif_{k}(\log(1+t))+(\log(1+t))^{m-1}\\
&\qquad \qquad \qquad \qquad \times (\partial_{t}(\log(1+t))Lif_{k}(\log(1+t)))\\
&=(\log(1+t))^{m-1}\frac{1}{1+t}\{(m-1)Lif_{k}(\log(1+t))+Lif_{k-1}(\log(1+t))\}.\\
\end{split}
\end{equation}

By (\ref{52}) and  (\ref{53}), we get

\begin{equation}\label{54}
\begin{split}
&\langle (\log(1+t))^{m}Lif_{k}(\log(1+t)) \vert x^{n} \rangle \\
&=(m-1)\langle Lif_{k}(\log(1+t))(1+t)^{-1} \vert (\log(1+t))^{m-1}x^{n-1} \rangle \\
&\qquad\qquad +\langle (Lif_{k-1}(\log(1+t))(1+t)^{-1} \vert (\log(1+t))^{m-1}x^{n-1} \rangle\\
&=\sum_{l=0}^{n-m}(m-1)!\binom{n-1}{l+m-1}S_{1}(l+m-1,m-1)\{(m-1)C_{n-l-m}^{(k)}(1)+C_{n-l-m}^{(k-1)}(1)\}.\\
\end{split}
\end{equation}

Therefore, by (\ref{51}) and  (\ref{54}), we obtain the following theorem.

\begin{theorem}\label{thm5}
For $n \geq m \geq 1$, we have
\begin{equation*}
\begin{split}
&\sum_{l=0}^{n-m}m!\binom{n}{l+m}S_{1}(l+m,m)C_{n-l-m}^{(k)}\\
&=\sum_{l=0}^{n-m}(m-1)!\binom{n-1}{l+m-1}S_{1}(l+m-1,m-1)\{(m-1)C_{n-l-m}^{(k)}(1)+C_{n-l-m}^{(k-1)}(1)\}.\\
\end{split}
\end{equation*}
\end{theorem}

In particular, for $n \geq 1$, we have

\begin{equation*}
C_{n-1}^{(k-1)}(1)=\sum_{l=0}^{n-1}(-1)^{l}l!\binom{n}{l+1}C_{n-l-1}^{(k)}.
\end{equation*}

From (\ref{1}), we note that

\begin{equation}\label{55}
\sum_{n=0}^{\infty}C_{n}\frac{t^{n}}{n!}=Lif_{1}(\log(1+t))=\frac{t}{\log(1+t)},
\end{equation}
where $C_{n}=C_{n}^{(1)}(0)$ is called the $n$-th Cauchy number of the first kind.

Let us consider the following sequences which are defined by the generating function to be

\begin{equation}\label{56}
\left(\frac{t}{\log(1+t)} \right)^{r}Lif_{k}(\log(1+t))=\sum_{n=0}^{\infty}T_{n}^{(r,k)}\frac{t^{n}}{n!}.
\end{equation}

Then, by (\ref{55}) and (\ref{56}), we get

\begin{equation}\label{57}
\left(\frac{t}{\log(1+t)} \right)^{r}Lif_{k}(\log(1+t))=\sum_{n=0}^{\infty}\left\{ \sum_{l_{1}+\cdots + l_{r+1}=n}\binom{n}{l_{1},\cdots , l_{r+1}}C_{l_{1}} \cdots C_{l_{r}}C_{l_{r+1}}^{(k)} \right\}\frac{t^{n}}{n!}.
\end{equation}

From (\ref{56}) and (\ref{57}), we have

\begin{equation}\label{58}
T_{n}^{(r,k)}=\sum_{l_{1}+\cdots + l_{r+1}=n}\binom{n}{l_{1},\cdots , l_{r+1}}C_{l_{1}} \cdots C_{l_{r}}C_{l_{r+1}}^{(k)}.
\end{equation}

For $n \geq 1$, by (\ref{10}), we get

\begin{equation}\label{59}
\begin{split}
C_{n}^{(k)}&=\langle Lif_{k}(\log(1+t) \vert x^{n} \rangle =\langle Lif_{k}(\log(1+t) \vert xx^{n-1} \rangle\\
&=\langle \partial_{t} (Lif_{k}(\log(1+t)) \vert x^{n-1} \rangle=\langle \frac{Lif_{k-1}(\log(1+t)-Lif_{k}(\log(1+t)}{(1+t)\log(1+t)} \vert x^{n-1} \rangle\\
&=\langle \frac{Lif_{k-1}(\log(1+t))-Lif_{k}(\log(1+t))}{(1+t)\log(1+t)} \vert \frac{1}{n}tx^{n} \rangle\\
&=\frac{1}{n}\langle \frac{t}{\log(1+t)}(Lif_{k-1}(\log(1+t))-Lif_{k}(\log(1+t))) \vert \frac{1}{1+t}x^{n} \rangle\\
\end{split}
\end{equation}

It is easy to show that

\begin{equation}\label{60}
\frac{1}{1+t}x^{n}=\sum_{l=0}^{\infty}(-t)^{l}x^{n}=\sum_{l=0}^{n}(-1)^{l}(n)_{l}x^{n-l}.
\end{equation}

By (\ref{59}) and (\ref{60}), we get

\begin{equation}\label{61}
\begin{split}
C_{n}^{(k)}&=\frac{1}{n}\sum_{l=0}^{n}(-1)^{l}(n)_{l}\langle \frac{t}{\log(1+t)}(Lif_{k-1}(\log(1+t))-Lif_{k}(\log(1+t))) \vert x^{n-l} \rangle\\
&=\frac{1}{n}\sum_{l=0}^{n}(-1)^{n-l}(n-l)!\binom{n}{l}(T_{l}^{(1,k-1)}-T_{l}^{(1,k)}).\\
\end{split}
\end{equation}

Therefore, by (\ref{61}), we obtain the following lemma.

\begin{lemma}\label{lem6}
For $k \in \mathbb{Z},~ n \geq  1$, we have
\begin{equation*}
C_{n}^{(k)}=\frac{1}{n}\sum_{l=0}^{n}(-1)^{n-l}(n-l)!\binom{n}{l}(T_{l}^{(1,k-1)}-T_{l}^{(1,k)}).
\end{equation*}
\end{lemma}

It is known that

\begin{equation}\label{62}
\begin{split}
&\partial_{t}^{m}Lif_{k}(\log(1+t))\\
&=\sum_{a=1}^{m}\sum_{l=0}^{a}S_{1}(m,a)S_{1}(a+1,l+1)\frac{Lif_{k-l}(\log(1+t))}{(1+t)^{m}(\log(1+t))^{a}},~~(see~[11]).
\end{split}
\end{equation}

For $n \geq m \geq 1$, by (\ref{62}), we get

\begin{equation}\label{63}
\begin{split}
C_{n}^{(k)}&=\langle Lif_{k}(\log(1+t)) \vert x^{n} \rangle=\langle \partial_{t}^{m}Lif_{k}(\log(1+t)) \vert x^{n-m} \rangle\\
&=\sum_{a=1}^{m}\sum_{l=0}^{a}S_{1}(m,a)S_{1}(a+1,l+1) \langle \frac{Lif_{k-l}(\log(1+t))}{(1+t)^{m}(\log(1+t))^{a}} \vert x^{n-m} \rangle\\
&=\sum_{a=1}^{m}\sum_{l=0}^{a}S_{1}(m,a)S_{1}(a+1,l+1) \frac{1}{(n-m+a)_{a}} \langle \frac{Lif_{k-l}(\log(1+t))}{(1+t)^{m}(\log(1+t))^{a}} \vert t^{a}x^{n-m+a} \rangle\\
&=\sum_{a=1}^{m}\sum_{l=0}^{a}S_{1}(m,a)S_{1}(a+1,l+1) \frac{1}{(n-m+a)_{a}}\\
 &\qquad\qquad\qquad \times \langle\left(\frac{t}{\log(1+t)}\right)^{a}Lif_{k-l}(\log(1+t)) \vert \frac{1}{(1+t)^{m}}x^{n-m+a} \rangle.\\
\end{split}
\end{equation}

Now, we observe that

\begin{equation}\label{64}
\frac{1}{(1+t)^{m}}x^{n-m+a}=\sum_{s=0}^{n-m+a}(-1)^{s}\binom{m+s-1}{s}(n-m+a)_{s}x^{n-m+a-s}.
\end{equation}

From (\ref{63}) and (\ref{64}), we have

\begin{equation}\label{65}
\begin{split}
C_{n}^{(k)}=\sum_{a=1}^{m}\sum_{l=0}^{a}\sum_{s=0}^{n-m+a}(-1)^{s}&\binom{m+s-1}{s}\frac{(n-m)!}{(n-m+a-s)!}\\
&\times S_{1}(m,a)S_{1}(a+1,l+1)T_{n-m+a-s}^{(a,k-l)}.\\
\end{split}
\end{equation}

Therefore, by (\ref{65}), we obtain the following lemma.

\begin{lemma}\label{lem7}
For $ n \geq m \geq 1$, we have
\begin{equation*}
\begin{split}
C_{n}^{(k)}=\sum_{a=1}^{m}\sum_{l=0}^{a}\sum_{s=0}^{n-m+a}(-1)^{s}&\binom{m+s-1}{s}\frac{(n-m)!}{(n-m+a-s)!}\\
&\times S_{1}(m,a)S_{1}(a+1,l+1)T_{n-m+a-s}^{(a,k-l)}.\\
\end{split}
\end{equation*}
\end{lemma}

For $S_{n}(x) \sim (g(t),f(t))$, it is known that

\begin{equation}\label{66}
\frac{d}{dx}S_{n}(x)=\sum_{l=0}^{n-1}\binom{n}{l} \langle \bar{f}(t) \vert x^{n-l} \rangle S_{l}(x).
\end{equation}

From (\ref{22}) and (\ref{66}), we have

\begin{equation}\label{67}
\begin{split}
\frac{d}{dx}C_{n}^{(k)}(x)&=\sum_{l=0}^{n-1}\binom{n}{l} \langle -\log(1+t) \vert x^{n-l} \rangle C_{l}^{(k)}(x)\\
&=\sum_{l=0}^{n-1}\binom{n}{l} \langle -\frac{\log(1+t)}{t}t \vert x^{n-l} \rangle C_{l}^{(k)}(x)\\
&=-\sum_{l=0}^{n-1}\binom{n}{l}(n-l) \langle \sum_{m=0}^{\infty}\frac{(-1)^{m}t^{m}}{m+1} \vert x^{n-l-1} \rangle C_{l}^{(k)}(x)\\
&=(-1)^{n}n!\sum_{l=0}^{n-1}\frac{(-1)^{l}}{(n-l)l!}C_{l}^{(k)}(x),~~(n \geq 1).
\end{split}
\end{equation}

For $C_{n}^{(k)}(x) \sim \left(\frac{1}{Lif_k(-t)},e^{-t}-1 \right),~B_{n}^{(r)}(x) \sim \left(\left(\frac{e^{t}-1}{t}\right)^{r},t\right)$, by (\ref{20}) and (\ref{21}), we have

\begin{equation}\label{68}
C_{n}^{(k)}(x)=\sum_{m=0}^{n}C_{n,m}B_{m}^{(r)}(x),
\end{equation}
where

\begin{equation}\label{69}
\begin{split}
C_{n,m}&=\frac{1}{m!} \langle \frac{\left(\frac{e^{-\log(1+t)}-1}{-\log(1+t)}\right)^{r}}{\frac{1}{Lif_k(\log(1+t))}}(-\log(1+t))^{m} \vert x^{n} \rangle \\
&=\frac{(-1)^{m}}{m!} \langle Lif_k(\log(1+t))\left(\frac{t}{(1+t)\log(1+t)}\right)^{r}(\log(1+t))^{m}\vert x^{n}\rangle \\
&=(-1)^{m}\sum_{l=0}^{n-m}\frac{1}{(l+m)!}S_{1}(l+m,m)(n)_{l+m} \\
&\qquad \qquad\qquad\times \langle Lif_k(\log(1+t))\left(\frac{t}{(1+t)\log(1+t)}\right)^{r} \vert x^{n-l-m}\rangle.\\
\end{split}
\end{equation}

Carlitz's polynomials $\beta_{n}^{(r)}(x)$ are defined by the generating function to be

\begin{equation}\label{70}
\left(\frac{t}{\log(1+t)}\right)^{r}(1+t)^{x}=\sum_{n=0}^{\infty}\beta_{n}^{(r)}(x) \frac{t^{n}}{n!},~~(see~[4,5]).
\end{equation}

By (\ref{69}) and (\ref{70}), we get

\begin{equation}\label{71}
C_{n,m}=(-1)^{m}\sum_{l=0}^{n-m}\sum_{a=0}^{n-m-l}\binom{n}{l+m}\binom{n-m-l}{a}S_{1}(l+m,m)\beta_{a}^{(r)}(-r)C_{n-m-l-a}^{(k)}.
\end{equation}

Therefore, by (\ref{68}) and (\ref{71}), we obtain the following theorem.

\begin{theorem}\label{thm8}
For $ k \in \mathbb{Z},~ n \geq 0$, we have
\begin{equation*}
\begin{split}
C_{n}^{(k)}(x)=&\sum_{m=0}^{n}\left\{(-1)^{m}\sum_{l=0}^{n-m}\sum_{a=0}^{n-m-l}\binom{n}{l+m}\binom{n-m-l}{a}\right.\\
& \left. \times S_{1}(l+m,m)\beta_{a}^{(r)}(-r)C_{n-m-l-a}^{(k)} \right\}B_{n}^{(r)}(x).\\
\end{split}
\end{equation*}
\end{theorem}

\begin{remark}
It is known that
\begin{equation}\label{72}
\frac{t}{(1+t)\log(1+t)}=\sum_{a=0}^{\infty}B_{a}^{(a)}\frac{t^{a}}{a!}.
\end{equation}
\end{remark}

Thus, by  (\ref{72}), we get

\begin{equation}\label{73}
\left(\frac{t}{(1+t)\log(1+t)}\right)^{r}=\sum_{a=0}^{\infty}\left(\sum_{a_{1}+ \cdots +a_{r}=a}\binom{a}{a_{1}, \cdots ,a_{r}}B_{a_{1}}^{(a_{1})} \cdots B_{a_{r}}^{(a_{r})}\right)\frac{t^{a}}{a!}.
\end{equation}

From (\ref{69}) and (\ref{73}), we can also derive

\begin{equation}\label{74}
\begin{split}
C_{n,m}&=(-1)^{m}\sum_{l=0}^{n-m}\sum_{a=0}^{n-m-l}\sum_{a_{1}+ \cdots +a_{r}=a}\binom{n}{l+m}\binom{n-m-l}{a}\binom{a}{a_{1}, \cdots ,a_{r}}\\
& \times S_{1}(l+m,m)B_{a_{1}}^{(a_{1})} \cdots B_{a_{r}}^{(a_{r})}C_{n-m-l-a}^{(k)}.\\
\end{split}
\end{equation}

By (\ref{68}) and (\ref{74}), we get

\begin{equation}\label{75}
\begin{split}
C_{n}^{(k)}(x)&=\sum_{m=0}^{n}\left\{(-1)^{m}\sum_{l=0}^{n-m}\sum_{a=0}^{n-m-l}\sum_{a_{1}+ \cdots +a_{r}=a}\binom{n}{l+m}\binom{n-m-l}{a}\binom{a}{a_{1}, \cdots ,a_{r}}\right.\\
& \left. \times S_{1}(l+m,m) B_{a_{1}}^{(a_{1})} \cdots B_{a_{r}}^{(a_{r})}C_{n-m-l-a}^{(k)}\right\}B_{n}^{(r)}(x).\\
\end{split}
\end{equation}

From (\ref{3}) and (\ref{22}), we consider the following two Sheffer sequences:

\begin{equation}\label{76}
C_{n}^{(k)}(x) \sim \left(\frac{1}{Lif_k(-t)},e^{-t}-1 \right),~~H_n ^{(r)} (x|\lambda) \sim \left(\left(\frac{e^t-\lambda}{1-\lambda}\right)^{r},t\right),
\end{equation}
where $r \in \mathbb{Z}_{\geq0}$.

Let us assume that

\begin{equation}\label{77}
C_{n}^{(k)}(x)=\sum_{m=0}^{n}C_{n,m}H_m ^{(r)} (x|\lambda).
\end{equation}

Then, by (\ref{21}), we get

\begin{equation}\label{78}
\begin{split}
C_{n,m}&=\frac{(-1)^{m}}{m!(1-\lambda)^{r}} \langle Lif_k(\log(1+t))\left(\frac{1}{1+t}-\lambda \right)^{r} \vert (\log(1+t))^{m}x^n \rangle\\
&=\frac{(-1)^{m}}{(1-\lambda)^{r}}\sum_{l=0}^{n-m}\binom{n}{l+m}S_{1}(l+m,m)\sum_{a=0}^{r}\binom{r}{a}(-\lambda)^{r-a}\\
&\qquad\qquad\qquad \times \langle Lif_k(\log(1+t))(1+t)^{-a} \vert x^{n-m-l} \rangle\\
&=\frac{(-1)^{m}}{(1-\lambda)^{r}}\sum_{l=0}^{n-m}\binom{n}{l+m}S_{1}(l+m,m)\sum_{a=0}^{r}\binom{r}{a}(-\lambda)^{r-a}C_{n-m-l}^{(k)}(a).\\
\end{split}
\end{equation}

Therefore, by (\ref{77}) and (\ref{78}), we obtain the following theorem.

\begin{theorem}\label{thm9}
For $ k \in \mathbb{Z}, ~n,r \in \mathbb{Z}_{\geq 0}$, we have
\begin{equation*}
\begin{split}
C_{n}^{(k)}(x)=\sum_{m=0}^{n}&\left\{\frac{(-1)^{m+r}}{(1-\lambda)^{r}}\sum_{l=0}^{n-m}\sum_{a=0}^{r}(-1)^{a}\binom{n}{l+m}\binom{r}{a}\right.\\
&\left. \times S_{1}(l+m,m)\lambda^{r-a}
C_{n-m-l}^{(k)}(a)\right\}H_{m}^{(r)}(x|\lambda).\\
\end{split}
\end{equation*}
\end{theorem}

\begin{remark}
By the same method, we can see that
\begin{equation}\label{79}
\begin{split}
C_{n}^{(k)}(x)&=\sum_{m=0}^{n}\left\{\frac{(-1)^{m+r}}{(1-\lambda)^{r}}\sum_{l=0}^{n-m}\sum_{b=0}^{n-m-l}\sum_{a=0}^{r}(-1)^{a+b}\binom{n}{l+m}\binom{r}{a}
\binom{a+b-1}{b}\right.\\
&\left. \times (n-m-l)_{b}\lambda^{r-a}S_{1}(l+m,m)C_{n-m-l-b}^{(k)}
\right\}H_{m}^{(r)}(x|\lambda).
\end{split}
\end{equation}
\end{remark}

For $C_{n}^{(k)}(x) \sim \left(\frac{1}{Lif_k(-t)},e^{-t}-1 \right)$, $x^{(n)} \sim (1, 1-e^{-t})$, we have

\begin{equation}\label{80}
C_{n}^{(k)}(x) =\sum_{m=0}^{n}C_{n,m}x^{(m)},
\end{equation}
where

\begin{equation}\label{81}
\begin{split}
C_{n,m}&=\frac{1}{m!} \langle Lif_k(\log(1+t))(-t)^{m} \vert x^n \rangle\\
&=\frac{(-1)^{m}}{m!}(n)_{m} \langle Lif_k(\log(1+t)) \vert x^{n-m} \rangle\\
&=(-1)^{m}\binom{n}{m}C_{n-m}^{(k)}.
\end{split}
\end{equation}

Therefore, by (\ref{77}) and (\ref{78}), we obtain the following theorem.

\begin{theorem}\label{thm10}
For $ k \in \mathbb{Z}, ~n\geq 0$, we have
\begin{equation*}
C_{n}^{(k)}(x)=\sum_{m=0}^{n}(-1)^{m}\binom{n}{m}C_{n-m}^{(k)}x^{(m)},
\end{equation*}
\end{theorem}
where $x^{(n)}=x(x+1)\cdots(x+n-1)$.


\bigskip
ACKNOWLEDGEMENTS. This work was supported by the National Research Foundation of Korea(NRF) grant funded by the Korea government(MOE)\\
(No.2012R1A1A2003786 ).
\bigskip


\end{document}